\input amstex
\magnification=1200
\documentstyle{amsppt}
\NoRunningHeads
\NoBlackBoxes
\topmatter
\title Kaleidoscope-roulette:\linebreak
the resonance phenomena in perception games
\endtitle
\author Denis V. Juriev\endauthor
\affil ul.Miklukho-Maklaya 20-180, Moscow 117437 Russia\linebreak
(e-mail: denis\@juriev.msk.ru)
\endaffil
\date math.OC/9905180
\enddate
\abstract\nofrills Kaleidoscope-roulettes, a proper class of perception games,
is described. Kaleidoscope-roulette is defined as a perception and, hence,
verbalizable interactive game, whose hidden dialogue consists of quasirandom 
sequences of ``words''. The resonance phenomena in such games and their
controlling are discussed.
\endabstract
\endtopmatter
\document
The mathematical formalism of interactive games, which extends one of ordinary
games [1] and is based on the concept of an interactive control, was recently
proposed by the author [2] to take into account the complex composition of 
controls of a real human person, which are often complicated couplings of 
his/her cognitive and known controls with the unknown subconscious behavioral 
reactions. This formalism is applicable also to the description of external 
unknown influences and, thus, is useful for problems in computer science 
(e.g. the semi-artificially controlled distribution of resources), 
mathematical economics (e.g. the financial games with unknown dynamical 
factors) and sociology (e.g. the collective decision making).

Recently, two proper classes of the interactive games were introduced:
the verbalizable interactive games [3] and the perception games [4].
The first class appeared as a result of the interactive game theoretical
description of dialogues as psycholinguistic phenomena and the second one
was obtained as such description of perception processes and the image
understanding. Each perception game is a verbalizable interactive game.

This article is devoted to a new proper subclass of interactive games,
namely, to the kaleidoscope-roulettes. Kaleidoscope-roulette is defined
as a perception and, hence, verbalizable interactive game, whose hidden
dialogue consists of quasirandom sequences of ``words''. The resonance
phenomena in such games are investigated.

Though kaleidoscope-roulettes are naturally associated with an
entertainment their real applications may be far beyond it due to their
origin.

\head I. Interactive games\endhead

\subhead 1.1. Interactive systems and intention fields\endsubhead

\definition{Definition 1 [2]} An {\it interactive system\/} (with $n$
{\it interactive controls\/}) is a control system with $n$ independent 
controls coupled with unknown or incompletely known feedbacks (the feedbacks
as well as their couplings with controls are of a so complicated nature that 
their can not be described completely). An {\it interactive game\/} is a game 
with interactive controls of each player.
\enddefinition

Below we shall consider only deterministic and differential interactive
systems. In this case the general interactive system may be written in the 
form:
$$\dot\varphi=\Phi(\varphi,u_1,u_2,\ldots,u_n),\tag1$$
where $\varphi$ characterizes the state of the system and $u_i$ are
the interactive controls:
$$u_i(t)=u_i(u_i^\circ(t),\left.[\varphi(\tau)]\right|_{\tau\leqslant t}),$$
i.e. the independent controls $u_i^\circ(t)$ coupled with the feedbacks on
$\left.[\varphi(\tau)]\right|_{\tau\leqslant t}$. One may suppose that the
feedbacks are integrodifferential on $t$.

\proclaim{Proposition [2]} Each interactive system (1) may be transformed
to the form (2) below (which is not, however, unique):
$$\dot\varphi=\tilde\Phi(\varphi,\xi),\tag2$$
where the magnitude $\xi$ (with infinite degrees of freedom as a rule) 
obeys the equation
$$\dot\xi=\Xi(\xi,\varphi,\tilde u_1,\tilde u_2,\ldots,\tilde u_n),\tag3$$
where $\tilde u_i$ are the interactive controls of the form $\tilde 
u_i(t)=\tilde u_i(u_i^\circ(t); \varphi(t),\xi(t))$ (here the dependence
of $\tilde u_i$ on $\xi(t)$ and $\varphi(t)$ is differential on $t$, i.e.
the feedbacks are precisely of the form 
$\tilde u_i(t)=\tilde u_i(u_i^\circ(t);\varphi(t),\xi(t),
\dot\varphi(t),\dot\xi(t),\ddot\varphi(t),\ddot\xi(t),\ldots,
\varphi^{(k)}(t),\mathbreak\xi^{(k)}(t))$).
\endproclaim

\remark{Remark 1} One may exclude $\varphi(t)$ from the feedbacks in
the interactive controls $\tilde u_i(t)$. One may also exclude the
derivatives of $\xi$ and $\varphi$ on $t$ from the feedbacks.
\endremark

\definition{Definition 2 [2]} The magnitude $\xi$ with its dynamical equations
(3) and its cont\-ri\-bution into the interactive controls $\tilde u_i$ will 
be called the {\it intention field}.
\enddefinition

Note that the theorem holds true for the interactive games. In practice, the 
intention fields may be often considered as a field-theoretic description of 
subconscious individual and collective behavioral reactions. However, they 
may be used also the accounting of unknown or incompletely known external 
influences. Therefore, such approach is applicable to problems of computer 
science (e.g. semi-automatically controlled resource distribution) or 
mathematical economics (e.g. financial games with unknown factors).
The interactive games with the differential dependence of feedbacks are
called {\it differential}. Thus, the theorem states a possibility of
a reduction of any interactive game to a differential interactive game
by introduction of additional parameters -- {\sl the intention fields}.

\subhead 1.2. Some generalizations\endsubhead The interactive games introduced 
above may be generalized in the following ways. 

The first way, which leads to the {\it indeterminate interactive games},
is based on the idea that the pure controls $u_i^\circ(t)$ and the 
interactive controls $u_i(t)$ should not be obligatory related in the
considered way. More generally one should only postulate that there are
some time-independent quantities $F_\alpha(u_i(t),u_i^\circ(t),\varphi(t),
\ldots,\varphi^{(k)}(t))$ for the independent magnitudes $u_i(t)$ and 
$u_i^\circ(t)$. Such claim is evidently weaker than one of Def.1. For 
instance, one may consider the inverse dependence of the pure and 
interactive controls: $u_i^\circ(t)=u_i^\circ(u_i(t),\varphi(t),\ldots,
\varphi^{(k)}(t))$.

The inverse dependence of the pure and interactive controls has a nice
psychological interpretation. Instead of thinking of our action consisting
of the conscious and unconscious parts and interpreting the least as 
unknown feedbacks, which ``dress'' the first, one is able to consider
our action as a single whole whereas the act of consciousness is in
the extraction of a part, which it declares as its property.

The second way, which leads to the {\it coalition interactive games}, is
based on the idea to consider the games with coalitions of actions and to
claim that the interactive controls belong to such coalitions. In this case
the evolution equations have the form
$$\dot\varphi=\Phi(\varphi,v_1,\ldots,v_m),$$
where $v_i$ is the interactive control of the $i$-th coalition. If the 
$i$-th coalition is defined by the subset $I_i$ of all players then
$$v_i=v_i(\varphi(t),\ldots,\varphi^{(k)}(t),u^\circ_j| j\in I_i).$$
Certainly, the intersections of different sets $I_i$ may be non-empty so
that any player may belong to several coalitions of actions. Def.1 gives the
particular case when $I_i=\{i\}$.

The coalition interactive games may be an effective tool for an analysis of
the collective decision making in the real coalition games that spread the
applicability of the elaborating interactive game theory to the diverse 
problems of sociology. 

\subhead 1.3. Differential interactive games and their 
$\varepsilon$--representations\endsubhead 

\definition{Definition 3 [3]} The {\it $\varepsilon$--representation\/} of 
differential interactive game is a representation of the differential
feedbacks in the form
$$u_i(t)=u_i(u_i^\circ,\varphi(t),\ldots,\varphi^{(k)}(t);
\varepsilon_i(t))\tag4$$
with the {\sl known\/} function $u_i$ of all its arguments, where
the magnitudes $\varepsilon_i(t)\in\Cal E$ are {\sl unknown\/} functions of
$u_i^\circ$ and $\varphi(t)$ with its higher derivatives:
$$\varepsilon_i(t)=\varepsilon_i(u_i^\circ(t),\varphi(t),\dot\varphi(t),
\ldots,\varphi^{(k)}(t)).$$
\enddefinition

It is interesting to consider several different $\varepsilon$-representations
simultaneously. For such simultaneous $\varepsilon$-representations
with $\varepsilon$-parameters $\varepsilon_i^{(\alpha)}$ a crucial role is
played by the time-independent relations between them:
$$F_\beta(\varepsilon_i^{(1)},\ldots,\varepsilon_i^{(\alpha)},\ldots,
\varepsilon_i^{(N)}; u_i^\circ,\varphi,\ldots,\varphi^{(k)})\equiv0,$$
which are called the {\it correlation integrals}. Certainly, in practice
the correlation integrals are determined {\sl a posteriori\/} and, thus they 
contain an important information on the interactive game. Using the sufficient
number of correlation integrals one is able to construct various algebraic 
structures in analogy to the correlation functions in statistical physics 
and quantum field theory.

\head II. Dialogues and the verbalizable interactive games.
Perception games.\endhead

\subhead 2.1. Dialogues as interactive games. The verbalization\endsubhead

Dialogues as psycholinguistic phenomena can be formalized in terms of
interactive games. First of all, note that one is able to consider
interactive games of discrete time as well as interactive games of
continuous time above.

\definition{Defintion 4A (the na{\"\i}ve definition of dialogues) [3]}
The {\it dialogue\/} is a 2-person interactive game of discrete time with 
intention fields of continuous time.
\enddefinition

The states and the controls of a dialogue correspond to the speech whereas 
the intention fields describe the understanding. 

Let us give the formal mathematical definition of dialogues now.

\definition{Definition 4B (the formal definition of dialogues) [3]}
The {\it dialogue\/} is a 2-person interactive game of discrete time of 
the form
$$\varphi_n=\Phi(\varphi_{n-1},\vec v_n,\xi(\tau)| 
t_{n-1}\!\leqslant\!\tau\!\leqslant\!t_n).\tag5$$
Here $\varphi_n\!=\!\varphi(t_n)$ are the states of the system at the
moments $t_n$ ($t_0\!<\!t_1\!<\!t_2\!<\!\ldots\!<\!t_n\!<\!\ldots$), 
$\vec v_n\!=\!\vec v(t_n)\!=\!(v_1(t_n),v_2(t_n))$ are the interactive 
controls at the same moments; $\xi(\tau)$ are the intention fields of 
continuous time with evolution equations
$$\dot\xi(t)=\Xi(\xi(t),\vec u(t)),\tag6$$
where $\vec u(t)=(u_1(t),u_2(t))$ are continuous interactive controls with 
$\varepsilon$--represented couplings of feedbacks:
$$u_i(t)=u_i(u_i^\circ(t),\xi(t);\varepsilon_i(t)).$$
The states $\varphi_n$ and the interactive controls $\vec v_n$ are certain
{\sl known\/} functions of the form
$$\aligned
\varphi_n=&\varphi_n(\vec\varepsilon(\tau),\xi(\tau)| 
t_{n-1}\!\leqslant\!\tau\!\leqslant\!t_n),\\
\vec v_n=&\vec v_n(\vec u^\circ(\tau),\xi(\tau)|
t_{n-1}\!\leqslant\!\tau\!\leqslant\!t_n).
\endaligned\tag7
$$
\enddefinition

Note that the most nontrivial part of mathematical formalization of dialogues
is the claim that the states of the dialogue (which describe a speech) are 
certain ``mean values'' of the $\varepsilon$--parameters of the intention
fields (which describe the understanding).

\remark{Important}
The definition of dialogue may be generalized on arbitrary number of players
and below we shall consider any number $n$ of them, e.g. $n=1$ or $n=3$, 
though it slightly contradicts to the common meaning of the word ``dialogue''.
\endremark

An embedding of dialogues into the interactive game theoretical picture
generates the reciprocal problem: how to interpret an arbitrary differential
interactive game as a dialogue. Such interpretation will be called the
{\it verbalization}.

\definition{Definition 5 [3]} A differential interactive game of the form
$$\dot\varphi(t)=\Phi(\varphi(t),\vec u(t))$$
with $\varepsilon$--represented couplings of feedbacks 
$$u_i(t)=u_i(u^\circ_i(t),\varphi(t),\dot\varphi(t),\ddot\varphi(t),\ldots
\varphi^{(k)}(t);\varepsilon_i(t))$$
is called {\it verbalizable\/} if there exist {\sl a posteriori\/}
partition $t_0\!<\!t_1\!<\!t_2\!<\!\ldots\!<\!t_n\!<\!\ldots$ and the 
integrodifferential functionals
$$\aligned
\omega_n&(\vec\varepsilon(\tau),\varphi(\tau)|
t_{n-1}\!\leqslant\!\tau\!\leqslant\!t_n),\\
\vec v_n&(\vec u^\circ(\tau),\varphi(\tau)|
t_{n-1}\!\leqslant\!\tau\!\leqslant\!t_n)
\endaligned\tag8$$ 
such that
$$\omega_n=\Omega(\omega_{n-1},v_n;\varphi(\tau)|
t_{n-1}\!\leqslant\!\tau\!\leqslant\!t_n).
\tag 9$$
\enddefinition

The verbalizable differential interactive games realize a dialogue in sense
of Def.4.

\remark{Remark 2} One may include $\omega_n$ explicitely into the evolution
equations for $\varphi$
$$\dot\varphi(\tau)=\Phi(\varphi(\tau),\vec u(\tau);\omega_n),\quad 
\tau\in[t_n,t_{n+1}].$$
as well as into the feedbacks and their couplings.
\endremark

The main heuristic hypothesis is that all differential interactive games
``which appear in practice'' are verbalizable. The verbalization means that 
the states of a differential interactive game are interpreted as intention 
fields of a hidden dialogue and the problem is to describe such dialogue 
completely. If a differential interactive game is verbalizable one 
is able to consider many linguistic (e.g. the formal grammar of a related 
hidden dialogue) or psycholinguistic (e.g. the dynamical correlation of 
various implications) aspects of it.

During the verbalization it is a problem to determine the moments $t_i$. A 
way to the solution lies in the structure of $\varepsilon$-representation.
Let the space $E$ of all admissible values of $\varepsilon$-parameters be
a CW-complex. Then $t_i$ are just the moments of transition of the 
$\varepsilon$-parameters to a new cell. 

\subhead 2.2. Perception games\endsubhead

\definition{Definition 6} The {\it perception game\/} is a multistage
verbalizable game (no matter finite or infinite) for which the intervals 
$[t_i,t_{i+1}]$ are just the sets. The conditions of their finishing 
depends only on the current value of $\varphi$ and the state of $\omega$ 
at the beginning of the set. The initial position of the set is the final 
position of the preceeding one.
\enddefinition

Practically, the definition describes the discrete character of the
perception and the image understanding. For example, the goal of a concrete
set may be to perceive or to understand certain detail of the whole image.
Another example is a continuous perception of the moving or changing object.

Note that the definition of perception games is applicable to various forms 
of perception, though the most interesting one is the visual perception.
The proposed definition allows to take into account the dialogical character
of the image understanding and to consider the visual perception, image
understanding and the verbal (and nonverbal) dialogues together. It may be
extremely useful for the analysis of collective perception, understanding 
and controlling processes in the dynamical environments -- sports, dancings, 
martial arts, the collective controlling of moving objects, etc. On the other 
hand this definition explicates the self-organizing features of human 
perception, which may be unraveled by the game theoretical analysis. And, 
finally, the definition put a basis for a systematical application of the 
linguistic (e.g. formal grammars) and psycholinguistic methods to the 
image understanding as a verbalizable interactive game with a mathematical 
rigor. 

\head III. Kaleidoscope-roulettes and the resonance phenomena\endhead

\subhead 3.1. Kaleidoscope-roulette\endsubhead
The kaleidoscope-roulette is a result of the attempt to combine the
kaleidoscope, one of the simplest and effective visual game, with the
roulette essentially using the elements of randomness and the treatment of
resonances. The main idea is to substitute random sequences of roulette
by the quasirandom sequences, which may be generated by the interactive
kaleidoscope. The obtained formal definition is below.

\definition{Definition 7} {\it Kaleidoscope-roulette\/} is a perception
game with a quasirandom sequence of quantites $\{\omega_n\}$.
\enddefinition

Certainly, the explicit form of functionals (8) is not known to the players.

Many concrete versions of kaleidoscope-roulettes are constructed. Though 
they are naturally associated with an entertainment their real applications 
may be far beyond it due to their origin and the abstract character of their 
definition. 

\subhead 3.2. The resonance phenomena in kaleidoscope-roulettes\endsubhead
Though the sequence $\{\omega_n\}$ is quasirandom the equations (9) for them
may have the resonance solutions. The resonance means a dynamical correlation
of two quasirandom sequences $\{v_n\}$ and $\{\omega_n\}$ whatever $\varphi$
is realized. In such case the quantities $\{v_n\}$ may be comprehended as
``fortune'', what is not senseless in contrast to the ordinary roulette.
However, $v_n$ are {\sl interactive\/} controls and their explicit 
dependence on $\vec u^\circ$ and $\varphi$ is not known. Nevertheless, one 
is able to use {\sl a posteriori\/} analysis and short-term predictions
based on it (cf.[5]) if the time interval $\Delta t$ in the short-term
predictions is not less than the interval $t_{n+1}-t_n$. To do it 
one should slightly improve constructions of [5] to take the discrete-time 
character of $v_n$ into account. It allows to perform the short-term 
controlling of the resonances in a kaleidoscope-roulette if they are observed.
The conditions of applicability of short-time predictions to the controlling
of resonances may be expressed in the following form: one should claim that
{\sl variations of the interactivity should be slower than the change of sets 
in the considered multistage game}. 

\remark{Remark 3} The possibility to control resonances by $v_n$ using its 
short-term predictions does not contradict to its quasirandomness, because
$v_n$ is quasirandom with respect to $v_{n-1}$ but not to $\varphi(\tau)$
($\tau\in[t_n,t_{n+1}]$).
\endremark

\head IV. Conclusions\endhead

Kaleidoscope-roulettes, a proper class of perception games, is described. 
They are defined as perception and, hence, verbalizable interactive games, 
whose hidden dialogue consists of quasirandom sequences of ``words''. 
The resonance phenomena in such games and their controlling are discussed.
A possibility of the short-term controlling of resonances in the 
kaleidoscope-roulettes is doubtless an intriguing feature for its use for 
the entertainment purposes as well as far beyond them.

\Refs
\roster
\item"[1]" Isaaks R., Differential games. Wiley, New York, 1965;\newline
Owen G., Game theory, Saunders, Philadelphia, 1968.
\item"[2]" Juriev D., Interactive games and representation theory. I,II.
E-prints: math.FA/9803020, math.RT/9808098.
\item"[3]" Juriev D., Interactive games, dialogues and the verbalization.
E-print: math.OC/9903001.
\item"[4]" Juriev D., Perception games, the image understanding and
interpretational geometries. E-print: math.OC/9905165.
\item"[5]" Juriev D., The laced interactive games and their {\sl a
posteriori\/} analysis. E-print:\linebreak math.OC/9901043; Differential
interactive games: The short-term predictions. E-print: math.OC/9901074.
\endroster
\endRefs
\enddocument